\documentclass[11pt,letterpaper]{article}

\usepackage{fullpage}

\usepackage{multicol}
\usepackage{wrapfig}

\usepackage[dvipsnames]{xcolor}
\usepackage{wallpaper}
\usepackage{graphicx}
\usepackage[margin=0.9cm, font=small, labelfont=sc]{caption}
\usepackage{subcaption}

\def\showauthornotes{0}

\def\showkeys{0}
\def\showdraftbox{0}
\def\showcolorlinks{1}

\def\showfixme{0}

\def\arxivmode{0}
\def\fastmode{0}

\ifnum\fastmode=0
\usepackage{etex}
\fi

\ifnum\fastmode=0
\usepackage[l2tabu, orthodox]{nag}
\fi

\usepackage{xspace,enumerate}

\ifnum\fastmode=0
\usepackage[T1]{fontenc}
\usepackage[full]{textcomp}
\fi

\usepackage[american]{babel}

\usepackage{mathtools}

\usepackage{amsthm}

\newtheorem{theorem}{Theorem}[section]
\newtheorem*{theorem*}{Theorem}

\newtheorem*{proposition*}{Proposition}

\newtheorem*{lemma*}{Lemma}

\newtheorem*{conjecture*}{Conjecture}

\newtheorem*{fact*}{Fact}

\newtheorem*{exercise*}{Exercise}

\newtheorem*{hypothesis*}{Hypothesis}

\theoremstyle{definition}

\newtheorem{exercise-easy}[theorem]{Exercise}
\newtheorem{exercise-med}[theorem]{Exercise}
\newtheorem{exercise-hard}[theorem]{Exercise$^\star$}

\newtheorem*{claim*}{Claim}

\newtheorem*{remark*}{Remark}

\newtheorem*{observation*}{Observation}

\ifnum\arxivmode=0
\usepackage{newpxtext} %
\usepackage{textcomp} %
\usepackage[varg,bigdelims]{newpxmath}
\usepackage[scr=rsfso]{mathalfa}%
\usepackage{bm} %
\let\mathbb\varmathbb
\fi

\ifnum\arxivmode=1
\usepackage[varg]{pxfonts} %
\usepackage{textcomp} %
\usepackage[scr=rsfso]{mathalfa}%
\usepackage{bm} %
\fi

\ifnum\showkeys=1
\usepackage[color]{showkeys}
\fi

\makeatletter
\@ifpackageloaded{hyperref}
{}
{\usepackage[pagebackref]{hyperref}}

\definecolor{bleudefrance}{rgb}{0.01, 0.1, 1.0}
\definecolor{azure}{rgb}{0.0, 0.5, 1.0}

\ifnum\showcolorlinks=1
\hypersetup{
pagebackref=true,
colorlinks=true,
urlcolor=blue,
linkcolor=blue,
citecolor=OliveGreen}
\fi

\ifnum\showcolorlinks=0
\hypersetup{
pagebackref=true,
colorlinks=false,
pdfborder={0 0 0}}
\fi

\usepackage{prettyref}

\newcommand{\savehyperref}[2]{\texorpdfstring{\hyperref[#1]{#2}}{#2}}

\newrefformat{eq}{\savehyperref{#1}{\textup{(\ref*{#1})}}}
\newrefformat{lem}{\savehyperref{#1}{Lemma~\ref*{#1}}}
\newrefformat{def}{\savehyperref{#1}{Definition~\ref*{#1}}}
\newrefformat{obs}{\savehyperref{#1}{Observation~\ref*{#1}}}
\newrefformat{ass}{\savehyperref{#1}{Assumption~\ref*{#1}}}
\newrefformat{thm}{\savehyperref{#1}{Theorem~\ref*{#1}}}
\newrefformat{cor}{\savehyperref{#1}{Corollary~\ref*{#1}}}
\newrefformat{cha}{\savehyperref{#1}{Chapter~\ref*{#1}}}
\newrefformat{sec}{\savehyperref{#1}{Section~\ref*{#1}}}
\newrefformat{app}{\savehyperref{#1}{Appendix~\ref*{#1}}}
\newrefformat{tab}{\savehyperref{#1}{Table~\ref*{#1}}}
\newrefformat{fig}{\savehyperref{#1}{Figure~\ref*{#1}}}
\newrefformat{hyp}{\savehyperref{#1}{Hypothesis~\ref*{#1}}}
\newrefformat{alg}{\savehyperref{#1}{Algorithm~\ref*{#1}}}
\newrefformat{rem}{\savehyperref{#1}{Remark~\ref*{#1}}}
\newrefformat{item}{\savehyperref{#1}{Item~\ref*{#1}}}
\newrefformat{step}{\savehyperref{#1}{step~\ref*{#1}}}
\newrefformat{conj}{\savehyperref{#1}{Conjecture~\ref*{#1}}}
\newrefformat{fact}{\savehyperref{#1}{Fact~\ref*{#1}}}
\newrefformat{prop}{\savehyperref{#1}{Proposition~\ref*{#1}}}
\newrefformat{prob}{\savehyperref{#1}{Problem~\ref*{#1}}}
\newrefformat{claim}{\savehyperref{#1}{Claim~\ref*{#1}}}
\newrefformat{relax}{\savehyperref{#1}{Relaxation~\ref*{#1}}}
\newrefformat{red}{\savehyperref{#1}{Reduction~\ref*{#1}}}
\newrefformat{part}{\savehyperref{#1}{Part~\ref*{#1}}}
\newrefformat{ex}{\savehyperref{#1}{Exercise~\ref*{#1}}}
\newrefformat{property}{\savehyperref{#1}{Property~\ref*{#1}}}

\newcommand{\Sref}[1]{\hyperref[#1]{\S\ref*{#1}}}

\usepackage{nicefrac}
\ifnum\showauthornotes=2
\newcommand{\mynotes}[1]{{\sffamily\small\color{teal}{#1}}\medskip}
\newcommand{\Authornote}[2]{{\sffamily\small\color{Maroon}{[#1: #2]}}\medskip}
\newcommand{\Authornotecolored}[3]{{\sffamily\small\color{#1}{[#2: #3]}}}
\newcommand{\Authorcomment}[2]{{\sffamily\small\color{gray}{[#1: #2]}}}
\newcommand{\Authorstartcomment}[1]{\sffamily\small\color{gray}[#1: }

\newcommand{\Authorfnote}[2]{\footnote{\color{red}{#1: #2}}}
\newcommand{\Authorfixme}[1]{\Authornote{#1}{\textbf{??}}}
\newcommand{\Authormarginmark}[1]{\marginpar{\textcolor{red}{\fbox{\Large #1:!}}}}
\newcommand{\myexplain}[1]{{\sffamily\small\color{red}{\noindent [Explanation:\medskip\newline \begin{quote}#1\hfill]\end{quote}}}\medskip}
\newcommand{\explain}[1]{{\sffamily\small\color{red}{#1}}\medskip}

\else

\newcommand{\mynotes}[1]{}
\newcommand{\Authornote}[2]{}
\newcommand{\Authornotecolored}[3]{}
\newcommand{\Authorcomment}[2]{}
\newcommand{\Authorstartcomment}[1]{}

\newcommand{\Authorfnote}[2]{}
\newcommand{\Authorfixme}[1]{}
\newcommand{\Authormarginmark}[1]{}
\newcommand{\myexplain}[1]{}
\newcommand{\explain}[1]{}

\ifnum\showauthornotes=1
\renewcommand{\myexplain}[1]{{\sffamily\small\color{red}{\noindent \begin{quote}{\bf Explanation:} \medskip\newline #1\end{quote}}}\medskip}
\fi

\fi

\ifnum\showfixme=0

\fi

\usepackage{boxedminipage}

\newcommand{\Esymb}{\mathbb{E}}

\newcommand{\Vsymb}{\mathbb{V}}

\DeclareMathOperator*{\E}{\Esymb}
\DeclareMathOperator*{\Var}{\Vsymb}

\newcommand{\textparen}[1]{\text{(#1)}}

\ifx\because\undefined
\newcommand{\because}[1]{\textparen{because #1}}
\else
\renewcommand{\because}[1]{\textparen{because #1}}
\fi

\newcommand\bdot\bullet

\ifx\mathds\undefined %

\else

\fi

\DeclareMathOperator{\argmax}{argmax}

\newcommand{\Z}{\mathbb Z}

\newcommand{\R}{\mathbb R}

\newcommand{\cA}{\mathcal A}
\newcommand{\cB}{\mathcal B}

\newcommand{\cL}{\mathcal L}

\newcommand{\cP}{\mathcal P}

\newcommand{\cS}{\mathcal S}

\newcommand{\cU}{\mathcal U}

\newcommand{\cW}{\mathcal W}

\newcommand{\bbH}{\mathbb H}

\newcommand{\bbZ}{\mathbb Z}

\newcommand{\bbP}{\mathbb P}

\renewcommand{\le}{\leqslant}

\renewcommand{\ge}{\geqslant}

\ifnum\showdraftbox=1

\else

\fi

\let\epsilon=\varepsilon

\numberwithin{equation}{section}

\newcommand\MYcurrentlabel{xxx}

\newcommand{\MYstore}[2]{%
  \global\expandafter \def \csname MYMEMORY #1 \endcsname{#2}%
}

\newcommand{\MYload}[1]{%
  \csname MYMEMORY #1 \endcsname%
}

\newcommand{\MYnewlabel}[1]{%
  \renewcommand\MYcurrentlabel{#1}%
  \MYoldlabel{#1}%
}

\newcommand{\MYdummylabel}[1]{}

\newcommand{\torestate}[1]{%
  \let\MYoldlabel\label%
  \let\label\MYnewlabel%
  #1%
  \MYstore{\MYcurrentlabel}{#1}%
  \let\label\MYoldlabel%
}

\newcommand{\restatetheorem}[1]{%
  \let\MYoldlabel\label
  \let\label\MYdummylabel
  \begin{theorem*}[Restatement of \prettyref{#1}]
    \MYload{#1}
  \end{theorem*}
  \let\label\MYoldlabel
}

\newcommand{\restatelemma}[1]{%
  \let\MYoldlabel\label
  \let\label\MYdummylabel
  \begin{lemma*}[Restatement of \prettyref{#1}]
    \MYload{#1}
  \end{lemma*}
  \let\label\MYoldlabel
}

\newcommand{\restateprop}[1]{%
  \let\MYoldlabel\label
  \let\label\MYdummylabel
  \begin{proposition*}[Restatement of \prettyref{#1}]
    \MYload{#1}
  \end{proposition*}
  \let\label\MYoldlabel
}

\newcommand{\restatefact}[1]{%
  \let\MYoldlabel\label
  \let\label\MYdummylabel
  \begin{fact*}[Restatement of \prettyref{#1}]
    \MYload{#1}
  \end{fact*}
  \let\label\MYoldlabel
}

\newcommand{\restate}[1]{%
  \let\MYoldlabel\label
  \let\label\MYdummylabel
  \MYload{#1}
  \let\label\MYoldlabel
}

\newcommand{\addreferencesection}{
  \phantomsection
\ifnum\stocmode=0
  \addcontentsline{toc}{section}{References}
\else
  \addcontentsline{toc}{section}{References \hspace*{1in} --------- End of extended abstract ---------}
\fi

}

\newcommand{\e}{\epsilon}

\renewcommand{\paragraph}[1]{\medskip\noindent{\bf #1.}}

\allowdisplaybreaks

\usepackage{paralist}

\renewcommand{\Vsymb}{\mathrm{Var}}

\newcommand{\vertiii}[1]{{\left\vert\kern-0.25ex\left\vert\kern-0.25ex\left\vert #1 
          \right\vert\kern-0.25ex\right\vert\kern-0.25ex\right\vert}}

\newcommand{\polytope}[1]{\text{\textup{\textsc{#1}}}}

\newcommand*{\corr}{\polytope{corr}}

\usepackage{enumitem}
\usepackage[normalem]{ulem}
\usepackage{stmaryrd}

\newcommand{\cmnt}[1]{}

\renewcommand{\mathbb}{\vvmathbb}

\date{}

\title{\vspace*{15mm}%
{Random metric geometries on the plane and Kardar-Parisi-Zhang universality}
\author{Shirshendu Ganguly \thanks{Associate Professor, Department of Statistics, University of California, Berkeley, \texttt{sganguly@berkeley.edu}}
}}

\begin{document}

\maketitle

\ThisULCornerWallPaper{1}{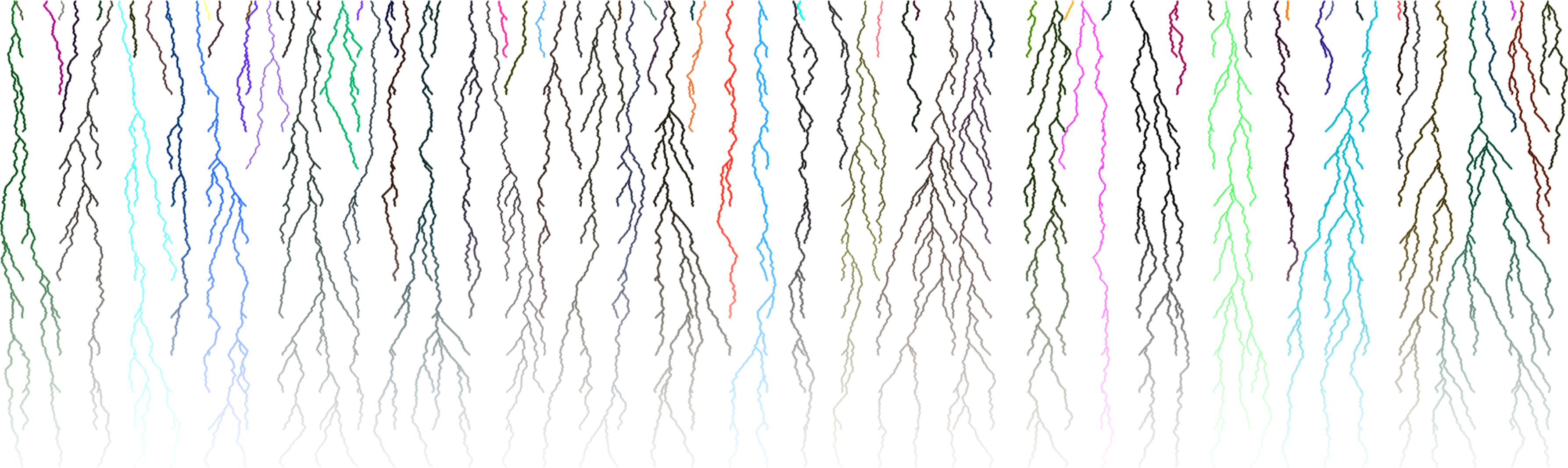}

\textit{A central object featuring in a significant aspect of modern research in probability theory is a random metric space obtained by distorting Euclidean space using some form of random noise. A topic of fundamental importance in this regard is then how the geodesics behave in the presence of noise.}

While examples of random metric spaces occur naturally in various physical situations, a few things have to be made precise before formulating concrete mathematical questions. For simplicity, we will work in a ``discrete'' setting, i.e., instead of $\R^d$, consider the lattice $\Z^d$ with the usual nearest neighbor graph metric. In the latter, the distance between two points is the smallest number of edges on any path connecting the same, i.e., this is simply the well known $\ell_1$ metric.  

We can now describe a canonical model of random geometry introduced in 1965 by Hammersley and Welsh.  The description is deceptively simple! For each edge $e\in \Z^d$, let $X_e$ be an i.i.d (independent and identically distributed) copy of a non-negative random variable with a fixed distribution, say $F.$ Thus, the ``edge weight'' $X_e$ 
can be thought of as the random length of the edge $e,$ which was previously simply one, and this leads to a new (random) shortest path metric on $\Z^d.$
Let $T(x,y)$ and $\Gamma(x,y)$ denote the random distance and the corresponding shortest path(s) between $x$ and $y$ respectively. 

\begin{figure}[h]
\centering
\includegraphics[width=.8\textwidth]{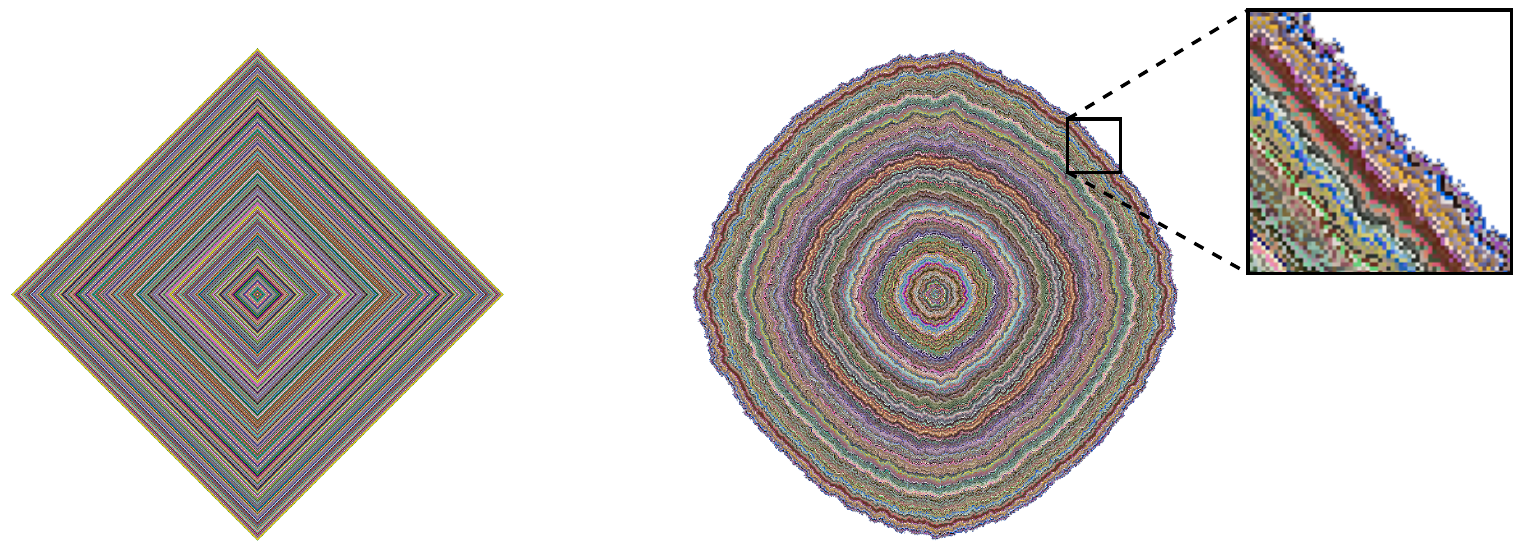}
\caption{Left: The $\ell_1$ ball around the origin in $\Z^2.$ Middle: The metric ball in FPP around the origin with noise formed by i.i.d. Exponential variables. Right: Zooming in on the boundary fluctuation behavior. The colors denote the level curves, i.e., points at the same distance from the origin.}\label{f.shape}
\end{figure}

Originally this was suggested as a model of fluid flow through a porous but random medium where for any edge $e$ its weight $X_e$ denotes the time taken by the fluid to cross it, often termed as the \emph{passage time}. Hence the model is called First Passage Percolation (FPP).  Subsequently this has been of interest in many areas of  physics, biology as well as computer science, modeling  various natural phenomena including spreading of bacteria and infection, propagation of forest fires and flame fronts, flow of current and more. 

It might be instructive to begin by listing a series of questions for FPP  that one might wish to answer (for several more, see \cite{FPPsurvey}).

\begin{enumerate} 
\itemsep-.2em
\item \textit{How does the metric ball centered at any point, say $x,$ i.e., $\cB_x(r)=\{y\in \Z^d, T(x,y)\le r\},$ grow}?
\item
\textit{ The above question can be broken down into two broad subparts. The first is about comparing the macroscopic behavior to the $\ell_1$ metric on $\Z^d$. How does $T(x,y)$ behave compared to $\|x-y\|_1?$  }
\item
\textit{ The next question is more refined and is about the order of fluctuations. How does $T(x,y)-\E(T(x,y))$ ($\E$ denotes the expectation)  behave? More generally, how does the boundary of the metric ball fluctuate? }
 \item 
\textit{ Moving on from $T(x,y),$ the next thing to investigate is $\Gamma(x,y),$ the geodesic(s) between $x$ and $y.$ What is its fluctuation? In particular, how much does it deviate from the straight line joining $x$ and $y$?}
 \item 
 \textit{Finally, is there any sense in which one can take the limit of the \textbf{entire} metric space? Is the limit universal, i.e., it does not depend on the edge weight distribution $F,$ as long as it is ``sufficiently nice''?}
 \end{enumerate}

Perhaps a bit surprisingly, it turns out that this 
simple mathematical model is notoriously hard to analyze with  most of the above questions in full generality still remaining open.  Nonetheless, impressive, though somewhat limited, advances have been made over the years giving birth to key new ideas and techniques in mathematics, statistical physics and probability theory. Recent years have also seen an explosion of activity around related models of random geometry admitting certain ``exact formulas''.

In this article we will review some of the general ideas, followed by the more recent developments around the specialized exactly solvable examples. 

\section{Limiting shape and fluctuations}\label{s:shape}
Starting with Questions 1 and 2 about the growth rate of the metric,  the classical law of large numbers result (LLN) states that for any i.i.d. sequence of random variables $X_1,X_2,\ldots$ with finite mean, i.e., $\E(|X_1|)< \infty,$ the long term average of the sequence converges to $\E(X_1),$ i.e.  ${\sum_{i=1}^n X_i}/{n}\to \E(X_1).$ This in particular implies that the sum of i.i.d. variables grows linearly if $\E(X_1)\neq 0$.

Given the above, it is natural to wonder if a counterpart result exists for FPP where the random metric is also driven by i.i.d. variables, albeit in a complicated way.  
This was answered by Richardson and further refined by Cox, Durrett and Kesten in the form of the following \emph{shape theorem} (see Figure \ref{f.shape}).   

\begin{theorem}\label{shape}(Informal) For any $F,$ {sufficiently nice},  there exists a deterministic compact convex set $\cS$  (which depends on $F$) such that for any $\e>0,$ 
$$\bbP\left((1-\e)\cS\subset\frac{\cB_0(t)}{t}\subset(1+\e)\cS \text{ for all large } t \right)= 1. $$
\end{theorem}

Thus, in words, the above says that the metric ball around the origin (or for that matter around any point) essentially looks like a dilated version (by a factor of $t$) of $\cS.$ In particular it says that for each unit  $\vec v \in \R^d,$ the passage time has an asymptotic speed $\mu_{\vec v}$ depending on the geometry of $\cS,$ i.e.,  
$\frac{T(0, n {\vec v})}{n}\to 1/\mu_{\vec v}.$ This answers the question about growth rate in full generality. Note that this a \emph{non-universal} result, as the limiting shape $\cS$ depends on the details of the edge weight distribution $F,$ although, perhaps surprisingly, the precise description is not explicitly known in any case. 

We now arrive at Question 3, about fluctuations, i.e., say for $e_1=(1,0,0, \ldots)$, how does $$T(0,ne_1)-\E(T(0,ne_1))$$ behave?
This seemingly innocuous question is one of the key problems in probability theory! 

To begin, let us start with a classical fluctuation counterpart of the above law of large numbers result.
Under the additional assumption $\E(X^2_1)< \infty,$ the well-known  central limit theorem (CLT) states that 
$$\frac{\sum_{i=1}^nX_i-n\E(X_1)}{\sqrt{n\Var(X_1)}} \to N(0,1),$$  
where $N(0,1)$ is a standard Gaussian random variable and the convergence is in distribution. Note that unlike the LLN result, the above is \emph{universal}, yielding the Gaussian distribution in the limit regardless of the distribution of the   variables $X_i.$

Beyond weak convergence, often in applications, it is useful to obtain concentration bounds of the following kind (which need further assumptions that we will not spell out).  There exists an universal constant $C>0$ such that for all $t\ge 0$  
\begin{equation}\label{tail}
\bbP\left(\left|\frac{\sum_{i=1}^nX_i-n\E(X_1)}{\sqrt{\Var(X_1)}}\right|\ge t\right)\le C\exp\left(-\frac{t^2}{Cn}\right).
\end{equation}

Concentration of measure for linear functions of i.i.d. variables as above is very well understood. However in many natural examples, such as FPP, the observables of interest are complicated non-linear functions of i.i.d. variables. 
Proving concentration bounds for such random variables poses a major challenge in probability theory with diverse applications. Over the years the theory has seen major advances.  
{In the context of FPP,  concentration results have been proven in landmark results of Kesten followed by Talagrand, further by Benjamini, Kalai, Schramm and then subsequently improved across various works (see e.g., \cite{FPPsurvey}}). Below we record a somewhat informal statement which captures the flavor of some of these results. 

\begin{theorem}(Informal and for simplicity stated only in the $e_1$ direction.)\label{fpptail} For FPP on $\Z^d,$ under certain {``tail''} conditions on the edge-weight distribution,  the following Gaussian tail bound holds:
\begin{align}\label{tail2}
\bbP\left(|T(0,ne_1)-\E(T(0,ne_1))| \ge t \sqrt{n} \right) \le 
C\exp(-t^2/C).\end{align} 
Subsequently, the following improved variance bound was established. $$\Var(T(0,ne_1))\le C\frac{n}{\log n}.$$
Above, $C$ is a universal constant.
 \end{theorem}
 Note that the above concentration result for $T(0,ne_1)$ at scale $\sqrt n$ is similar to the case of the sum of i.i.d. variables in \eqref{tail}. 
However, it is easy to verify that the minimum of a family of i.i.d. random variables is more concentrated than each of the individual variables. Thus we should expect $T(0,ne_1)$, being the minimum over the lengths of all paths between $0$ and $ne_1$, to be more concentrated than the length of any individual path which is  a sum of i.i.d. random variables!  

Towards this, observe that while the Gaussian tail bound immediately implies an order $n$ bound on the variance, the second result above has a logarithmic improvement. Though perhaps not apparent, even this seemingly small improvement involves completely novel breakthrough ideas around the notion of \emph{hyper-contractivity}, and \emph{is the best known bound in any generality currently!} More recently, corresponding improvements of \eqref{tail2} at scale $\sqrt{{n}/{\log n}}$ have been established.

Nonetheless, using non-rigorous methods one can in fact predict that $T(0,ne_1)$ is concentrated at scale $n^{1/3}$ instead of $n^{1/2}$. Verifying the same mathematically is a driving force behind much of the advance in this area of probability. Most of the discussion in the remainder of the article will be focussed on this and related issues, eventually addressing Questions 4 and 5 as well.

\section{Kardar-Parisi-Zhang universality}As indicated in the discussion following Theorem \ref{shape}, while the growth rate is non-universal in the shape theorem, the fluctuation theory is expected to be universal and  one might wonder what the counterpart object for FPP is, analogous to the Gaussian distribution appearing in the CLT.  
It turns out, in dimension 2 (recall we have been considering FPP on $\Z^d$ for a general $d$ so far) there is a very precise description of the predicted fluctuation behavior. Thus we will subsequently focus on the special case of the plane. 
The general predictions are expected to hold for a large class of examples modeling random growth beyond metric balls in random metric spaces. 

To begin, we first introduce another model of stochastic growth which is closely related to FPP,  but with some key differences, and with a similar name, Last Passage Percolation (LPP). Although expected to behave similarly, it turns out that the latter exhibits some extremely useful and surprising properties absent in FPP, which will be apparent soon.

\begin{figure}[h]
\hspace{1.5cm}
	\begin{subfigure}[t]{0.30\textwidth}
		\centering
		\includegraphics[width=\textwidth]{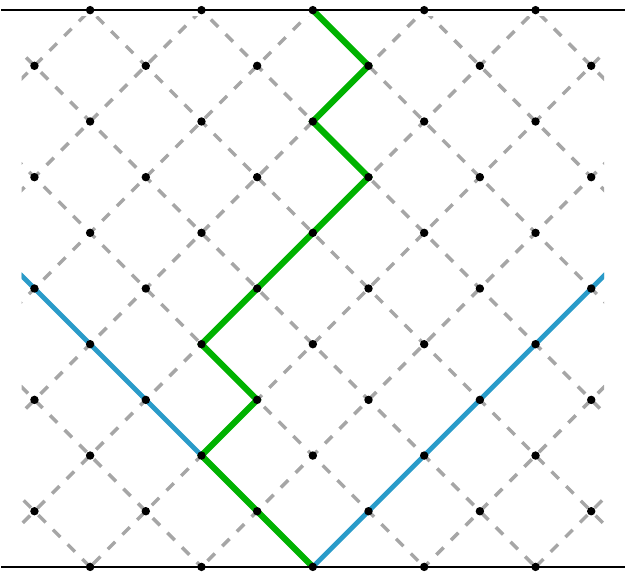}
	\end{subfigure} 
	\hspace{2.1cm}
	\begin{subfigure}[t]{0.43\textwidth}
		\centering
		 \raisebox{6mm}{\includegraphics[width=\textwidth]{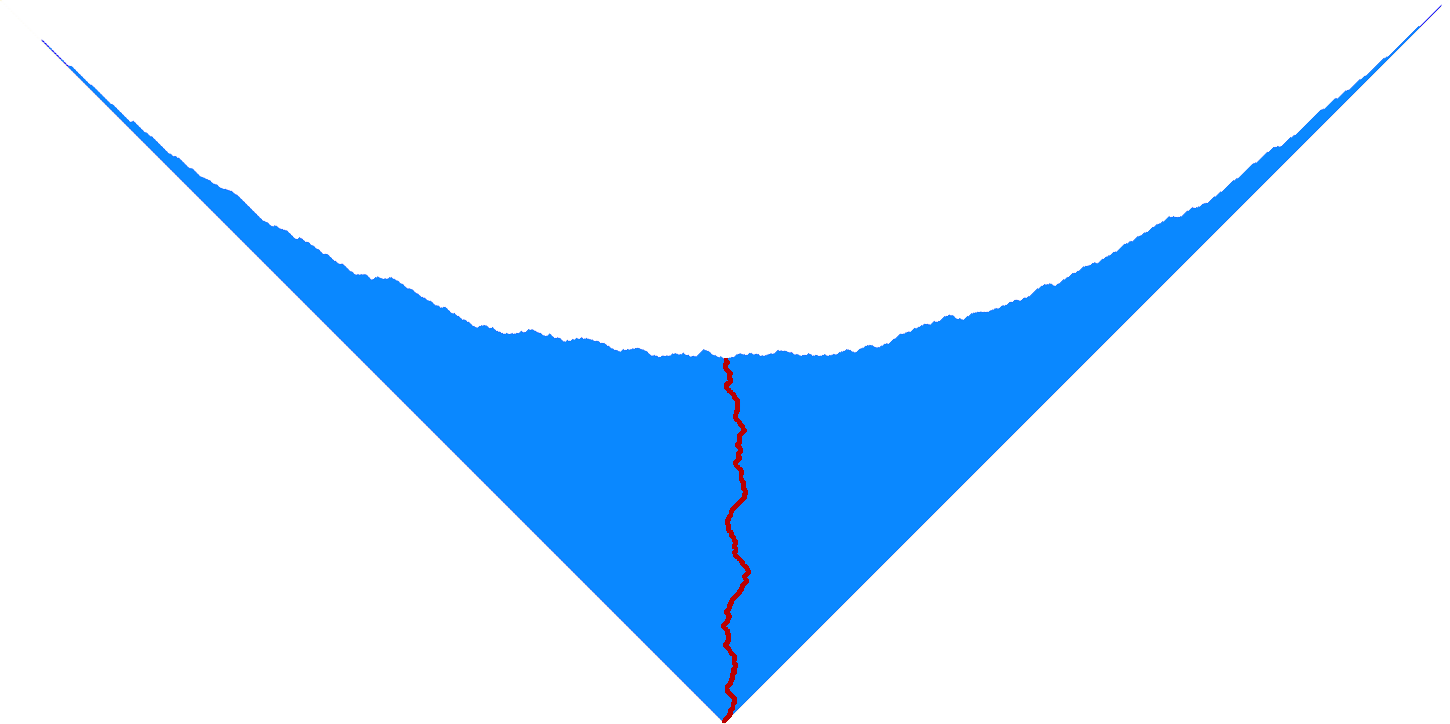}}
	\end{subfigure}
	\caption{Left: Last passage percolation on the upper half plane $\bbH^2_{+}$, with the horizontal direction as ``space" and vertical direction as ``time". The cone emanating from the origin denotes the set of points which can be connected to the origin using oriented paths. Right: Growth cluster ) of the origin, where the boundary, properly scaled converges to a ``truncated'' {hyperbola}. A geodesic between the origin and a point on the boundary of the cluster is illustrated.}\label{f.LPP}
\end{figure}

\textbf{Last Passage Percolation:} Consider the lattice upper half plane $\bbH^2_{+}=\{\big((x-y)/2^{1/2},(x+y)/2^{1/2}\big): (x,y) \in \bbZ^2, x+y \ge 0\}.$ formed by rotating $\Z^2$ by $45^\circ$ as shown in Figure \ref{f.LPP} (it will be convenient later to work with this orientation). Equip each vertex $v$ (instead of edges as in FPP) with a non-negative i.i.d. variable $X_v$ with a common distribution, say $F.$
Further, unlike FPP, given two points $x,y \in \bbH^2_+,$ we will only consider \emph{oriented} paths from $x$ to $y$ which only move north-east or north-west as in Figure \ref{f.LPP}. For any such  path $\gamma,$ associate to it the weight  $L(\gamma)$ by summing up the random variables $X_v$ along the same. Finally, instead of considering minimum passage times as in FPP, we consider maximum ones, i.e., define the last passage time $L(x,y)$ as $\sup_{\gamma: x\to y} L(\gamma).$ 

This slight alteration of the setting, i.e., considering oriented paths and maximum passage times, along with certain special choices of the vertex weights, will lead to quite remarkable algebraic properties which form the central underpinning of most of the significant advances in our understanding of such examples. Going forward,  abusing terminology, we will call the maximizing path(s)  the geodesic(s) between $x,y.$ 
Note that owing to length maximization, the LPP weights $L(\cdot,\cdot)$ don't form a metric but rather an anti-metric, i.e., it satisfies the reverse triangle inequality, $L(x,y)+L(y,z)\le L(x,z)$.
The reader is encouraged to review the concept of directed metric introduced by Dauvergne and Vir\'ag \cite{DV2021} which unifies the notions of metric and anti-metric.

As in FPP, one considers the growing cluster around any point, say the origin (to be denoted throughout by $0$), $\cU_r(0)=\{x\in \bbH^2_+: L(0,x)\le r\}.$ The arguments of the shape theorem for FPP go through to show the existence of a limiting function $\cS: \R\to \R$ such that with high probability the boundary of $\cU_{r}(0)$ is sandwiched between the graphs of $\cS$ scaled by $r(1\pm \e)$ respectively. 

On simulating either FPP or LPP on the computer, one observes that the boundary of the growing cluster around the origin exhibits certain key features (it will be useful to view 2 dimensions as one ``spatial" direction and one ``time" direction, see Figure \ref{f.LPP}). This includes a global smoothening phenomenon involving faster growth in rougher portions on the boundary than smoother parts. Further, crucially, the growth rate is expected to depend non-linearly on the local gradient.  Finally, as is obvious from the descriptions of the models, the local fluctuations is driven by i.i.d. noise. 

In 1986, physicists  Kardar, Parisi and Zhang (KPZ) \cite{KPZ86} put forward a unified fluctuation theory predicting that the behavior of stochastic growth exhibiting the above features, e.g., the boundary of the growth cluster in LPP, should be the same as that of a canonical non-linear stochastic PDE (the KPZ equation). For more on the KPZ equation, see e.g., \cite{quastel,corwin1}.

While several models are expected to be in the KPZ universality class, keeping with our theme of random metric spaces and geometry of geodesics, we will simply be focussing on FPP and LPP throughout the article, as they enjoy special geometric properties missing in the other examples.

\subsection{Characteristic exponents governing fluctuation.} \label{corrstructure}
Much of the predictions about the fluctuation theory for such models can be summarized in the triple of ``critical'' exponents $(1/3:2/3:1)$. 
To describe this, consider the LPP value $L(0,(0,n))$. This grows at linear speed, as evident from the shape theorem, and hence is of order $n$ (explaining the exponent $1$). However as emphasized before, the most interesting aspect is its fluctuation!

Answering this and explaining the $1/3$ exponent, it is predicted that 
$$\big|L(0,(0,n))-\E(L(0,(0,n)))\big|\approx n^{1/3},$$  
confirming the intuition that this is much more concentrated than the weight of any given path joining the points, which has $O(n^{1/2})$ fluctuations. 

The exponent $2/3$ is related to both correlation and  geodesic behavior. 
Firstly, $L(0,(x_1,n))$ and $L(0,(x_2,n))$ exhibit non-trivial correlation as long as $|x_1-x_2|\approx n^{2/3}.$
\begin{wrapfigure}[15]{}{4 cm}
\centering
\includegraphics[width=.18\textwidth]{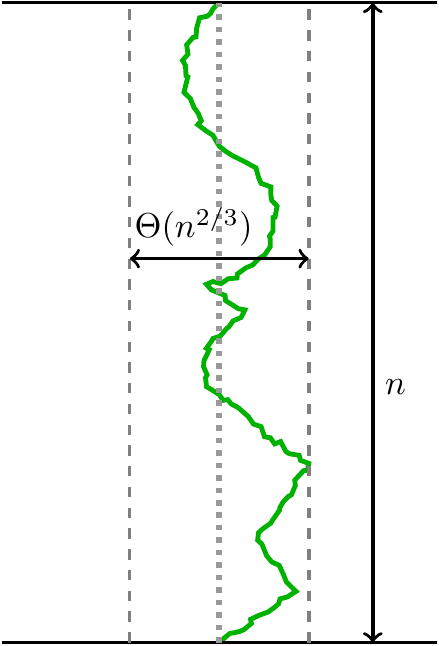} 
\caption{Illustrating the transversal fluctuation of $\Theta(n^{2/3})$ of the geodesic connecting $0$ and $(0,n).$}
\label{fig:TF}
\end{wrapfigure}
In particular, $L(0,(x_1,n))$ and $L(0,(x_2,n))$ are asymptotically (as $n\to \infty$) same if $|x_1-x_2|\ll n^{2/3}$ and independent if $|x_1-x_2|\gg n^{2/3}.$ 
Moreover, the exponent $2/3$ also shows up rather 
nicely in terms of geodesic behavior answering Question 4 stated at the beginning, about its transversal fluctuations (see Figure \ref{fig:TF}). Namely, the geodesic $\gamma$ between say $0$ and $(0,n)$ fluctuates by $\Theta(n^{2/3})$ away from the straight line joining its endpoints. As the reader might already be aware, a ``uniformly'' chosen path between the same endpoints  fluctuates by $\Theta(n^{1/2})$.
{Thus the geodesic has a bigger fluctuation than a random path, which at a very informal level is
explained by the fact that it is inclined to fluctuate more on its quest to find vertices of high weights.}

However, despite these intriguing predictions being around for multiple decades, the rigorous literature is lagging behind significantly as we have seen in FPP, where the best known bound on fluctuations is still $O(\sqrt n)$ up to a nontrivial logarithmic improvement, a far cry from the conjectured $O(n^{1/3})$ behavior!

Nonetheless, while the general situation is admittedly somewhat disappointing, a handful of examples possess certain special properties that have opened the door to various external mathematical tools leading to spectacular progress over the last twenty years. This brings us into the world of integrable probability. 

\section{Integrable Probability}

A class of models remarkably exhibit additional algebraic structure, though often quite difficult to discern, which can be exploited to obtain exact expressions for the observables of interest. The formulas can then be analyzed to extract information about their asymptotic behavior.  These include connections to representation theory, algebraic combinatorics and in particular to permutations, random matrices and their eigenvalues, quantum integrable systems and so on. 
To convey the basic idea, we will discuss a particular example of how a certain special choice of the vertex weights in LPP makes it \emph{integrable} or \emph{exactly solvable}. For more examples, see \cite{borodin2016lectures}.

\noindent
$\bullet$ \textbf{Exponential LPP:}\,\,  This is a model of LPP when the vertex weights $X_v$ are i.i.d. Exponential variables, i.e., with distribution $\bbP(X_v \ge y)=e^{-y}$ for $y\ge 0.$
 
To describe a key consequence of the algebraic properties of this model, we need to recall a  classical random matrix ensemble.

\noindent
$\bullet$ \textbf{Laguerre Unitary Ensemble (LUE):} For any positive integers $n,m,$ consider the $n\times m$ matrix $X$ with each entry being complex and i.i.d. with $X_{ij}=X^1_{ij}+\sqrt{-1} X^2_{ij}$
where $X^1_{ij},X^2_{ij},$ are i.i.d. standard Gaussian variables of variance $1/2$. The  LUE with parameters  $(m,n)$ is then given by $W=X^*X.$

We now discuss Johansson breakthrough observation: Exponential LPP is related to LUE and this can be used to verify KPZ predictions!

\subsection{$1/3$ exponent guiding weight fluctuations via random matrix theory}
The following remarkable identity for Exponential LPP was established by Johansson \cite{johansson1}: 
\begin{equation}
L\big(0, 2^{1/2}(m-n,m+n)\big)\overset{d}{=}\lambda_1
\end{equation}
 where $\lambda_1$ is the largest eigenvalue of the LUE with parameters $m,n$. He then analyzed the latter relying on the theory of orthogonal polynomials to prove the following seminal result. 

\begin{equation}\label{onept}
\frac{L(0,(0, 2^{1/2}n))-4n}{2^{4/3}n^{1/3}}\overset{d}{\to} F_{GUE}
\end{equation}
where $F_{GUE}$ is the well known Tracy-Widom distribution. This  describes the limiting fluctuation behavior of the largest eigenvalue of the Gaussian unitary ensemble (GUE), which analogous to LUE is another classical Hermitian random matrix  whose entries are, up to conjugation, i.i.d. complex Gaussians. Interestingly, $F_{GUE}$ has the following non-Gaussian tail behavior: as $\lambda \to \infty,$ $$F_{GUE}\big((\lambda,\infty)\big)\approx e^{-(4/3+o(1))\lambda^{3/2}},\,\,F_{GUE}\big((-\infty,-\lambda)\big)\approx e^{-(1/12+o(1)) \lambda^{3}}.$$

Thus this surprising turn of events connecting random planar growth and random matrices not only rigorously establishes the predicted $1/3$ exponent for the weight fluctuation for Exponential LPP, but also tells us that $F_{GUE}$ is the universal counterpart of the Gaussian distribution to be expected in this context, which no one has found a way to predict by other means. 

The first such result was proven by Baik, Deift and Johansson \cite{BDJ} in the context of a different model of LPP which is naturally connected to the well known problem of the longest increasing subsequence in a random permutation formulated by Ulam. We point the reader to the beautiful book \cite{romik} discussing in depth the above problem.

Having established the $1/3$ weight fluctuation exponent, the next natural goal is to obtain multipoint correlation information with the goal of verifying the $2/3$ exponent. This will also allow us to predict how the boundary of the metric ball in FPP behaves, addressing Question 3 listed at the beginning. 
Again, we will see the power and usefulness of the surprising random matrix connections!   

\subsection{$2/3$ exponent and spatial correlation structure}\label{s:corrstr}
Continuing from Section \ref{corrstructure}, to understand the correlation structure of $L(0,\cdot),$ \emph{at scale $n^{2/3}$,} it will be convenient to consider the scaled \emph{geodesic weight profile}, 
\begin{equation}\label{weightprofile}
\cL_n(x):=2^{-4/3}n^{-1/3}[L(0, 2^{1/2}(x2^{5/3}n^{2/3},n))-4n],
\end{equation}
 which encodes the properly centered and scaled last passage times from $0$ as the second endpoint is varied along the line $y=n$, generalizing the one point expression appearing in \eqref{onept}.

Again, it might be instructive to draw analogy with the classical case of a sequence of i.i.d. variables. Towards this, we recall Donsker's invariance principle, a process version of the CLT result. Namely, if one considers the random function $\cW_{n}(\cdot):[0,1]\to \R$ obtained by setting $$\cW_n(0)=0, \,\,\cW_n(j/n)={\sum_{i=1}^{j}[X_i-\E(X_i)]}/\sqrt{n \Var(X_1)}\,\, \text{for}\,\, 1\le j\le n,$$ and linearly interpolating for other points, then the process $\cW_n$ converges to a one-dimensional Brownian motion. The latter is a random process with independent increments, whose one point distributions are Gaussians. It is not an overstatement to claim that Brownian motion is the most universally occurring process in probability theory and is often the first guess for what the scaling limit of a random process in nature should be!

\begin{figure}[h]
\centering
\includegraphics[width=0.95\textwidth]{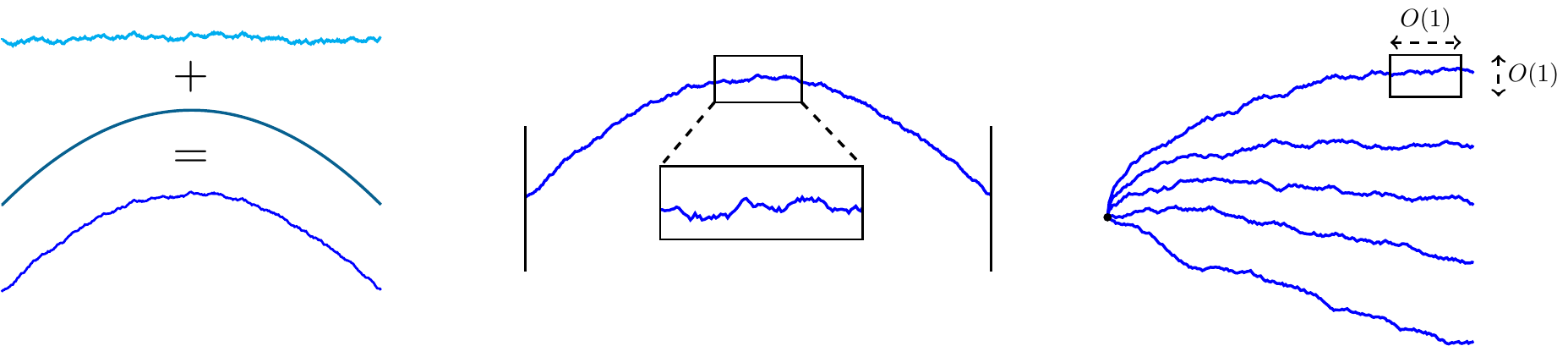}
\caption{Left: The stationary Airy$_2$ process minus a parabola gives $\cP\cA,$ the Parabolic Airy$_2$ process.  Middle: Zooming in on a compact interval, the latter shows Brownian behavior. Right: Illustrates Dyson Brownian motion, whose suitable edge scaling limit yields $\cP\cA$.}\label{f.airy line}
\end{figure}

However, we have already seen from \eqref{onept} that $\cL_n(0)$ converges to $F_{GUE}$, which is not Gaussian, and hence the process $\cL_n(\cdot)$ cannot converge to Brownian motion. 
In fact, it was shown using analysis of exact formulas (see e.g., \cite{borodin2008large, johansson2003discrete}), that the finite dimensional distributions of $\cL_n$ converge to that of a new process which has been termed as the \emph{Parabolic Airy$_2$ process}, which we will denote as (see Figure \ref{f.airy line}) $$\cP\cA:=\cA_{2}(x)-x^2.$$  
This is the central object in the KPZ universality class expected to  arise as the universal scaling limit of fluctuations in a large class of stochastic planar growth models including random metric spaces! 

Here, $\cA_2(x)$, called the Airy$_2$ process (without the parabolic correction), is a stationary ergodic process, whose one point distribution is the Tracy-Widom distribution. A key feature of $\cA_2$ is that it is \emph{determinantal} which admits usage of algebraic tools.
Without dwelling much on this, we simply mention that the latter means that the joint distribution of the values taken at different points is given by determinants of certain matrices. This in particular tells us how the correlation of the geodesic weight behaves at scale $n^{2/3},$ establishing mathematically the exponent $2/3$ in the triple $(1/3:2/3:1).$
 
The Airy$_2$ process was first constructed in the context of poly-nuclear growth, by Pr\"ahofer and Spohn \cite{prahofer2002scale}  who analyzed its spatial correlations using methods from the statistical mechanics of {Fermionic systems}. 

Finally, despite all the above acting as evidence that the first naive guess of every scaling limit being Brownian motion does not hold in this case, as we will now see, it wasn't too far!

\subsection{Local Brownianity of the Airy$_2$ process}Though the Airy$_2$ process is not literally Brownian motion, the fact that the LPP models are driven by i.i.d. variables indicates that the Airy$_2$ process ought to exhibit certain Brownian characteristics.  In fact this is a theme at the forefront of current research with major applications. In particular, this allows one to transfer our geometric knowledge of Brownian motion developed over a century to the Airy$_2$ process. 

The first result in this direction was proven by H\"agg \cite{hagg}. Subsequently, Corwin and Hammond \cite{CH14}  showed that for the Parabolic Airy$_2$ process $\mathcal{PA}$, a single sample of the increment process on any compact interval $[c,d],$ i.e., the process $\cP\cA(x+c)-\cP\cA(c)$ on the interval $[0,d-c]$ is \emph{indistinguishable} from a single sample of a properly scaled Brownian motion (see Figure \ref{f.airy line}).
However one cannot hope to have such an indistinguishability result on the whole real line since Airy$_2$ is a stationary process while Brownian motion is diffusive, i.e., it is approximately $\sqrt{x}$ at location $x$. 

The key tool in the proof? Again remarkable connections to random matrices.
The argument in the Corwin-Hammond work relied on a gorgeous identity proved by {O'Connell-Yor} \cite{o2002} which states that for Brownian LPP (a variant of Exponential LPP, where instead of sums of Exponential variables one has Brownian motions) the geodesic weight profile, also known to converge to $\cP\cA,$ admits an embedding as the top line in an ensemble of $n$ lines, such that the joint law of the $n$ lines is what is known as Dyson Brownian motion \cite{dyson}, or $n$ Brownian motions conditioned to avoid each other! This is a central object in random matrix theory, denoting the evolution of the eigenvalues of an $n\times n$ GUE whose entries evolve with time as independent Brownian motions.

As an immediate application, since the same is true for Brownian motion, one concludes that the Airy$_2$ process is almost surely nowhere differentiable and in fact is H\"older $1/2-$ continuous! A fact that could be rather complicated to establish via different means.

Thus, for certain models of LPP we have seen that the geodesic weight profile, scaled using KPZ exponents, converges to the Parabolic Airy$_2$ process which is locally Brownian like. This information can be used to analyze the boundary of the growth clusters in the corresponding models. 

As the reader might have noticed, in all the results outlined so far, connections to random matrices play a major role, and these highly non-obvious features only exist because of certain special choices made while defining the models, such as using Exponential variables or Brownian motions. There are other cases which are not directly connected to matrices but to other objects such as permutations but in the interest of brevity we will not be reviewing them and instead point the reader to \cite{borodin2016lectures}. 

We now come to Question 5 about the scaling limit of the metric space itself. Towards this, so far we have seen that when one point is fixed to be the origin, the distance to another point as the latter varies along a straight horizontal line, scales to $\cP\cA$. One is naturally led to wonder if one can take a similar scaling limit of the entire metric space itself, when one freely varies both the endpoints. A recent breakthrough addresses this.  

\subsection{The Directed Landscape} 
To talk about such scaling limits, we need to first upgrade our notation $\cL_n$ used in \eqref{weightprofile}, i.e., for $x,y\in \R$ and $0<s<t$, let $$\cL_n(y,s,x,t):=2^{-4/3}n^{-1/3}\left[L\left(2^{1/2}\big(y2^{5/3}n^{2/3}, \lfloor ns\rfloor\big),2^{1/2}\big(x2^{5/3}n^{2/3},\lfloor nt\rfloor\big)\right)-4n(t-s)\right],$$ which encodes the scaled geodesic weight between scaled points $(y,s)$ and $(x,t)$. (Note that the $4n$ term in \eqref{weightprofile} got replaced by $4n(t-s)$). As mentioned in Figure \ref{f.LPP}, we will view the vertical coordinates $s,t$ as time and the horizontal coordinates $x,y$ as space.

In 2018, Dauvergne, Ortmann and Virag \cite{DOV18} constructed the  Directed Landscape: a four parameter random energy field ${\cL(y,s,x,t): \R^4\to \R},$ which the pre-limiting field $\cL_n(y,s,x,t)$ is expected to converge to. The convergence to the same for $\cL_n$ and their counterparts in other integrable models have been proved across  \cite{DOV18} and subsequent work by Dauvergne and Vir\'ag \cite{DV2021}, making it conjecturally the universal scaling limit of such models of random geometry in two dimensions. 
Further, importantly, it was also shown that the geodesics in the pre-limiting models converge to their limiting counterparts in the Directed Landscape.

Having obtained the Directed Landscape $\cL$ as the scaling limit of the anti-metric structure in LPP, 
which is a rich universal object exhibiting intricate geometric features, unearthing the latter is in fact a topic of massive current interest. 
Before reviewing some of the recent advances in this direction, we begin by recording a few of its key properties:
\begin{itemize}
\itemsep-.2em
\item \textit{Translation invariance}: $\cL(y,s;x,t){=}\cL(0,0;x-y,t-s)$ in law. 
\item
\textit{Invariance under KPZ $(1/3:2/3:1)$ scaling}:  $\cL(y,s,x,t){=}q^{-1/3}\cL(q^{2/3}y,qs,q^{2/3}x,qt)$ in law, for $q>0$. \item \textit{Independence of increments}: for any time points $s_1<t_1<s_2<t_2\ldots <s_k<t_k,$ the processes $\cL(\cdot, s_i,\cdot, t_i)$ are independent across $i$.  
\end{itemize}
Note that $\cL$ contains several embeddings of the process $\cP\cA$. For instance, for any $y,$ $ \cL(y,0,y+\cdot,1){=}\cP\cA(\cdot)$ in law. This should be interpreted as the scaled geodesic weight profile which starts from $y$ instead of $0$. In fact the random two dimensional field $\cA(y,x):=\cL(y,0,x,1)$ is termed as the Parabolic-Airy sheet, providing a coupling of all the Parabolic-Airy$_2$ processes rooted at the different starting points. A lot has been recently proved  about this coupling structure.

\section{Geodesic geometry and its consequences}
Recall that the article started by indicating that one of its main motivations was to study geodesics in random distortions of the Euclidean metric.  All of the previous preparation now allows us to dive into geodesic geometry in integrable models of LPP, an area that has seen an explosion of activity recently.  We start by recording some basic but fundamental properties.

\noindent
$\bullet$ \textbf{Transversal fluctuations.}  
Recall that the transversal fluctuation of the geodesic from the straight line was mentioned to be governed by the exponent $2/3$. In fact, more can be said in terms of the process $\cP\cA.$ To see this consider Exponential LPP and the geodesic $\Gamma_{2n}$ from $0$ to $(0,2n)$ and say we are interested in the typical location of the point $v_*$ where $\Gamma$ intersects the line $y=n.$ Now such a $v_*$ maximizes $T\big(0,v\big)+T\big(v,(0,2n)\big)$ over all $v\in \{y=n\}.$ Now using independence of the noise field, it follows that, properly scaled, $T(0,\cdot)$ and $T(\cdot,(0,2n))$ converge to two independent $\cP\cA$ processes. Thus $v_*$ typically behaves as $(1+o(1))n^{2/3} x_*$ where $$x_*=\argmax_x \cP\cA^{(1)}(x)+\cP\cA^{(2)}(x)$$ and $\cP\cA^{(1)},\cP\cA^{(2)}$ are two independent Parabolic Airy$_2$ processes. 
(That transversal fluctuations are dictated by the KPZ exponent of $2/3$ was first rigorously established by Johansson \cite{johansson2}.)

\noindent
$\bullet$ \textbf{Coalescence of geodesics.} We now come to perhaps the most striking feature of the random metric spaces being discussed, standing sharply in contrast to Euclidean geometry. This is the phenomenon of \emph{coalescence}. As the name suggests, this means that geodesics between pairs of points that are not very far away from each other tend to meet and share a non-trivial amount of their journey (see e.g., the figure on the first page which illustrates a network of geodesics!). This is very different from Euclidean geometry  where geodesics are straight lines and meet at a single point or not at all.
In very brief, the reason for the above is that geodesics are weight maximizing paths, and hence each path, regardless of its endpoints, prefer to pass through the vertices whose random weights are the highest, leading  them to meet each other, i.e., coalesce.

While the initial progress in understanding properties of random anti-metrics arising in LPP, such as their fluctuations, was mostly through algebraic methods, 
the past few years have seen several advances with a focus on geodesic behavior, based primarily on geometric and probabilistic analysis. While it is impossible to review all of them, we include a sample list below to provide the reader a flavor. 
\subsection{Applications}
\noindent
$\bullet$ \textbf{Fractal geometry.}
Fractals are self-similar sets; they look the same at all scales. 
\begin{wrapfigure}{r}[-4pt]{4 cm}
\centering
\includegraphics[width=.20\textwidth]{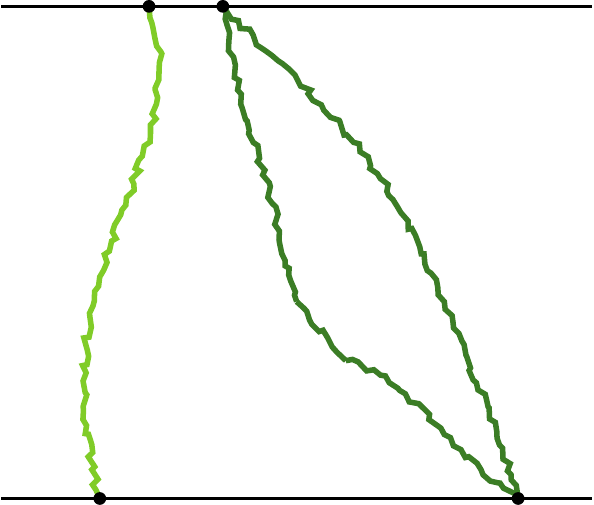} 
\caption{Unique or multiple disjoint geodesics between points.}
\label{fig:busemann}
\end{wrapfigure}
They are ubiquitous in nature and in mathematics, with the famous Cantor set being a classical example.
 Random fractal sets are also commonly occurring in probability. 
A good way is to think of them as rather sparse random subsets of $\R^d$, say, whose distribution as you zoom in or zoom out stays the same. Given such a set, a natural mathematical problem is to quantify how sparse it is. This is where the notion of the ``dimension'' of such a set comes in. Without going into precise definitions, informally it measures how many balls of size $\e$ (for some small $\e$) in the underlying metric space (which will be $\R^d$ for us) is needed to cover it. E.g.,  the unit cube $[0,1]^d\subset \R^d$ requires roughly $(1/\e)^{d}$  balls of size $\e$ indicating that its dimension is $d.$

The scale invariance satisfied by the Directed Landscape makes it a rich source of intricate fractal behavior. 
We now describe a particular example related to the behavior of geodesics. It turns out that for any two fixed points $(y,0)$ and $(x,1)$, almost surely, the geodesic between the two points in $\cL$ is unique. 
Nonetheless, on account of there being uncountably many possible endpoints, a natural set to consider is the set of points $(x,y)$ such that there are two \emph{disjoint} geodesics between the points $(y,0)$ and $(x,1).$ A sequence of  recent articles (see e.g., \cite{BGH19, BGH20, GH21} and the references therein) study this and related sets and in particular shows that this random set has dimension $1/2$ almost surely! More recently, there have been results comparing such sets to a canonical $1/2$ dimensional set, namely the set of zeros  of a one dimensional Brownian motion.\\

\noindent
$\bullet$ \textbf{Temporal correlation.} It is of much interest to understand how values at different points in space and time in the Directed Landscape are correlated. While determinantal formulas have provided a reasonable understanding of the spatial correlation observables such as $\corr(\cL(0,0,0,1),\cL(0,0,x,1)),$  i.e., the geodesic weights with one end fixed and the other end moving in the spatial direction, formulas for similar quantities in the temporal direction, say  $\corr(\cL(0,0,0,1),\cL(0,0,0,t))$ are much less amenable to asymptotic analysis. 
Predictions about the behavior of the latter as a function of $t$ were made in \cite{FS16}. 
These have a clear interpretation in terms of geodesic coalescence and we outline below the simplest case of $\corr(\cL(0,0,0,1),\cL(0,0,0,t))$ for short times, i.e. when $t\ll 1$. Consider the geodesics $\Gamma_s$ going from $0$ to $(0,s)$. When $t\ll 1,$ the geodesic $\Gamma_1$ overlaps with a macroscopic fraction of $\Gamma_t.$ Since the weight of $\Gamma_1$ beyond height $t$ depends on independent noise (recall the independence of increments of the directed landscape), for heuristic purposes one can assume that the latter contributes negligibly to the covariance. Thus the entire covariance comes from the overlap and hence is expected to have the same order as the variance of $\cL(0,0,0,t)$. 
Arguments of this type have been made precise recently (see e.g. \cite{BG18, BGZ20} and the references therein). \\

\noindent
$\bullet$ \textbf{Non-existence of bi-geodesics.} As we have seen, geodesics behave very differently in these types of random metric spaces compared to Euclidean spaces. Another related notion in which they are expected to differ pertains to the existence of bi-geodesics, i.e. a bi-infinite path such that every finite segment of it is a geodesic. Note that any straight line is a bi-geodesic in Euclidean space, whereas a central problem in FPP on $\Z^2,$ first posed by Furstenberg, is to show that almost surely there are no bi-geodesics. The question has gained fame through its connection to the
existence of non-trivial ground states of the two-dimensional Ising ferromagnet with random
exchange constants. To see, roughly, why bi-geodesics are implausible, let us simply consider the possibility that a bi-geodesic passes through the origin $(0,0).$ The latter implies that, for any large $n,$ there exists points $x_n$ and $y_n$ such that $|x_n|, |y_n| \approx n$ with the geodesic $\Gamma(x_n,y_n)$ passing through the origin. Now because of transversal fluctuation of the latter being expected to be of order $n^{2/3},$ the intersection of the geodesic with the line $y=0$ should be roughly uniformly distributed on an interval of size $n^{2/3}$ and hence it containing the origin should have probability approximately $n^{-2/3}.$ Since $n$ can be arbitrarily large, this shows the probability that a bi-geodesic passing through $0$ has probability $0.$ To make this precise, one has to observe that, by coalescence, effectively one has only finitely many such choices for $x_n,y_n$. Recently the above strategy, suggested by Newman, was implemented rigorously for Exponential LPP by Basu-Hoffman-Sly \cite{BHS2018}. (A separate argument was also provided by Bal\'azs, Busani and  Sepp\"al\"ainen \cite{BBS20} using different methods). It must also be pointed out that despite the lack of a comparable understanding of geodesics, the problem for FPP has also seen some impressive progress. We point the reader to \cite{FPPsurvey} for more on this. \\

\noindent
$\bullet$ \textbf{Large deviations.} While we have discussed various \emph{typical} features of random metric spaces, such as fluctuations of the metric and geodesics, guided by KPZ exponents and so on, there is still a non-trivial probability they behave very strangely. This is in fact an important topic in probability, i.e., to study the occurrence  of rare events or, as it is commonly called, large deviations.
In the context of LPP, for instance, this amounts to considering when the geodesic weight is atypically large or small. E.g., for Exponential LPP,  since typically $T(0,(0,n))=4n+O(n^{1/3})$ (recall from \eqref{onept}) one studies the events  $T(0,(0,n)) \ge 5n$ (upper tail) or  $\le 3n$ (lower tail) (moderate deviation events such as $T(0,(0,n))\ge 4n+ n^{1/2}$ have also been studied). These probabilities were very precisely computed first by Sepp\"al\"ainen using connections to particle systems and subsequently by Johansson using bijections to eigenvalues. However, often in the study of large deviations, it is of interest to understand how the system behaves conditioned on being atypical. It was shown recently that in contrast to the typical behavior of $O(n^{2/3}),$ the transversal fluctuation of the corresponding geodesic becomes $O(n^{1/2})$  (localized) for the upper tail by Basu-Ganguly \cite{BG19}, and $O(n)$ (delocalized) for the lower tail by Basu-Ganguly-Sly \cite{BGS17A}, respectively. 

At a very intuitive level, the discrepancy between the two behaviors is rooted in the fact that the upper tail simply demands one single path of high weight, while the lower tail forces \emph{all} paths to have low weight. \\

\noindent
$\bullet$ \textbf{Chaos.}  The final topic we touch upon relates to how a random system behaves as the underlying randomness gets perturbed over time. In many complex disordered systems this leads the phenomenon of \emph{chaos}. Without providing formal definitions, let us just mention that 
this means that several models in statistical 
mechanics are characterized by intricate energy
 landscapes 
 where the ground state, the configuration 
 with the lowest energy, is expected to  
alter profoundly when the disorder of the
model is slightly perturbed (see \cite{superconcentration}).

In the context of FPP or LPP, very informally this amounts to saying that while given any fixed endpoints there is typically a unique geodesic between them, there are many ``near'' geodesics whose weights are very close to that of the geodesic. Now a slight perturbation of the underlying noise variables causes the geodesic to jump from its current position to one of these competing candidates.

There have been conjectures about how models in KPZ should react to dynamical perturbation in the physics literature, but it was only recently that the understanding of a particular model of LPP was rich enough to be able to undertake a mathematical investigation of its dynamical aspects. In \cite{GHchaos}, a phase transition was established which, very informally, asserts that  perturbation smaller than a certain threshold does not alter the geodesic significantly while beyond that the geodesic exhibits very little similarity to its initial state.

\section{Major research directions.}
Having reviewed some of the recent and not so recent progress in understanding properties of random distortions of the Euclidean lattices using probabilistic and geometric ideas and often relying on crucial integrable inputs, we end with a brief discussion on some major research directions in this area. 
\begin{itemize}
\item Fluctuation bounds: While it seems rather difficult to verify the KPZ predictions without any algebraic structure, a major goal is to significantly improve the current state of the art for FPP. While there have been several impressive results under certain unproven assumptions (see in particular work of Chatterjee \cite{chatterjeeKPZ} verifying the so called `KPZ-relation' between the weight and transversal fluctuation exponents under the assumption that they exist in some strong sense), the set of completely unconditional results is somewhat limited.  Let us just state the following representative question in this direction pertaining to the variance bounds for the passage times in FPP: 

{\centering\textit{Show that for FPP in $\Z^2,$ under ``reasonable" assumptions, $n^{\e}\le \Var(T(0,ne_1))\le n^{1-\e}$ for some $\e>0.$}}

 It is worth pointing out that while the general known upper bound as stated earlier is $O({n}/{\log n}),$ the best known lower bound is only $O(\log n)$!

\item High dimensions:  While the shape theorem and concentration results from Section \ref{s:shape} hold in any dimensions, the KPZ exponent triple $(1/3:2/3:1)$ prediction holds only in $2=1+1$ (one space and one time) dimensions. Thus, while there has been significant progress on the plane, at least for the integrable examples, the situation for higher dimensions is much less advanced with lack of consensus about what to expect even in the non-rigorous literature. In the three dimensional case, certain special growth models (anisotropic growth) admit projections to two dimensional growth which were analyzed in work of Borodin and Ferrari \cite{anisotropic}. More recently, the KPZ equation in $2+1$ dimensions was analyzed in certain regimes of the parameter space (see e.g., \cite{KPZ2}). 
\item Non-i.i.d disorder:  One of the key features of a model that lead to KPZ universal behavior is that growth is driven by noise that de-correlates fast in space and time and so throughout the article we have considered examples where the underlying disorder distorting the Euclidean geometry is independent in space-time. While certain results known for FPP with i.i.d. disorder carry over to the setting where the noise is simply ergodic and translation invariant, the fluctuation theory is expected to heavily depend on the correlation structure of the noise. There has recently been an explosion of activity in the planar case (see e.g. \cite{LQGsurvey}), when the noise comes from the exponential of the Gaussian free field, a canonical two dimensional log-correlated Gaussian process, which can be considered as the natural higher dimensional analogue of Brownian motion. 
\end{itemize}

\textbf{Summary}: While this article attempted to review the recent and past advances in understanding random planar metric spaces with a focus on geodesic geometry, there are several other aspects of this theory, and stochastic growth in general, that were overlooked. We strongly encourage the reader to refer to the excellent surveys \cite{FPPsurvey, corwin1, quastel} and the references therein to gain a more comprehensive insight about this fascinating area of modern research. 

\section{Acknowledgements}The author thanks two anonymous referees for many detailed comments which significantly improved the article. The author also thanks Milind Hegde, Ella Hiesmayr, Adam Jaffe, Kyeongsik Nam and Lingfu Zhang for helpful feedback on earlier drafts. Special thanks goes to Milind for all the help, with the figures and much more, throughout the writing of this article.  The author was supported by NSF grant DMS-1855688, NSF Career grant DMS-1945172, and a Sloan Fellowship.

\bibliographystyle{alpha}
{\footnotesize
\bibliography{notices}}

\end{document}